\input amssym.def
\input amssym
\magnification=1200
\parindent0pt
\hsize=16 true cm
\baselineskip=13  pt plus .2pt
$ $

\def\Z{\Bbb Z}

\def\A{\Bbb A}

\centerline {\bf  ${\rm SL}(n,\Bbb Z)$ cannot act on small spheres}

\bigskip \bigskip

\centerline {Bruno P. Zimmermann}

\bigskip

\centerline {Universit\`a degli Studi di Trieste}
\centerline {Dipartimento di
Matematica e Informatica}
\centerline {34100 Trieste, Italy}
\centerline
{zimmer@units.it}

\vskip 1cm

Abstract. {\sl  The group ${\rm SL}(n,\Bbb Z)$ admits a smooth faithful action on
$S^{n-1}$, induced from its linear action on $\Bbb R^n$. We show that, if $m <
n-1$ and $n> 2$, any smooth action of ${\rm SL}(n,\Bbb Z)$ on a
mod 2 homology $m$-sphere, and in particular on the sphere $S^m$, is trivial.}

\bigskip \bigskip

{\bf 1. Introduction}

\medskip

By [We] any smooth action of ${\rm SL}(n,\Bbb Z)$ on the torus $T^m$ is trivial
if $m < n$ (whereas it acts linearly and faithfully on $T^n = \Bbb R^n/\Bbb
Z^n$). Also, by [FM] any continuous action of ${\rm SL}(n,\Bbb Z)$ on a closed
surface is trivial for sufficiently large values of $n$ (depending on the genus
of the surface). Our main result is the following analogue for smooth actions of
${\rm SL}(n,\Bbb Z)$ on mod 2 homology spheres:

\bigskip

{\bf Theorem 1}. {\sl  Let $n > 2$ and $m < n-1$. Any smooth$^1$ action of
${\rm SL}(n,\Bbb Z)$ on a mod 2 homology $m$-sphere, and in particular on
$S^m$, is trivial.}

\bigskip

We note that the group ${\rm SL}(n,\Bbb Z)$ admits a smooth faithful action on
$S^{n-1}$ (induced from its linear action on $\Bbb R^n$). Theorem 1 was
conjectured by Parwani [P] who proved that, if $m < n-1$ and $n>2$, any
smooth$^1$ action of ${\rm SL}(n,\Bbb Z)$ on a mod 2 homology $m$-sphere
factors through the action of a finite group. The main point of the proof in
[P] is an application of the Margulis finiteness theorem which implies that,
for $n>2$, ${\rm SL}(n,\Bbb Z)$ is almost simple that is every normal subgroup
is either finite and central, or of finite index. Considering the subgroup
$(\Bbb Z_2)^{n-1}$  of diagonal matrices with all diagonal entries equal to
$\pm 1$, Smith fixed point theory implies that some noncentral element of ${\rm
SL}(n,\Bbb Z)$ has to act trivially on the mod 2 homology $m$-sphere if $m <
n-1$, and hence by Margulis' theorem the kernel of the action has finite index.
By the result in [P], Theorem 1 is then a consequence of the following result
to be proved in section 2.

\medskip

---------------------------------------------------

$^1$ See the note at the end of the paper

\vfill  \eject

{\bf Theorem 2}. {\sl  Let $n > 2$ and $m < n-1$.  Any smooth action of a
finite quotient of ${\rm SL}(n,\Bbb Z)$ on a mod 2 homology $m$-sphere is trivial.}

\bigskip

We note that Theorem 1 is proved in [P] for actions of ${\rm
SL}(n,\Bbb Z)$ on $S^1$ and $S^2$. For the case of $S^1$, Witte [Wi] has shown
that every continuous action of a subgroup of finite index in ${\rm SL}(n,\Bbb
Z)$ on the circle $S^1$ factors through a finite group action. In the context
of the Zimmer program for actions of irreducible lattices in semisimple Lie
groups  of $\Bbb R$-rank at least two ([Z]), it has been conjectured by Farb
and Shalen [FS] that any smooth action of a finite-index subgroup of ${\rm
SL}(n,\Bbb Z)$, $n > 2$, on a compact $m$-manifold factors through the action
of a finite group if $m < n-1$; however this remains still open for actions
e.g. on spheres (since the proof in [P] relies heavily on the existence of
certain elements of finite order in ${\rm SL}(n,\Bbb Z)$ which has torsion-free
subgroups of finite index, it does not apply to this more general situation).

\bigskip

In general, it is not easy to decide whether a given finite group admits a
faithful action on a homology $m$-sphere, for smaller values of $m$. For
example, Theorem 2 implies that the simple groups of type ${\rm
PSL}(n,\Z/p\Z)$, $p$ prime, do not act faithfully on a mod 2 homology
$m$-sphere whenever $m<n-1$ and $n>2$; however, the following remains open:

\bigskip

Question. What is the minimal dimension of an integer (a mod 2) homology sphere
which admits a faithful action of

i) ${\rm PSL}(n,\Z/p\Z)$?

ii) the alternating group $\A_n$?

\bigskip

In the context of Theorem 1,  the group ${\rm SL}(n,\Bbb Z)$ has a finite
subgroup isomorphic to the alternating group $\A_{n+1}$, and it is easy to see
that, if the subgroup $\A_{n+1}$ does not act faithfully, then the whole group
${\rm SL}(n,\Bbb Z)$ has to act trivially (cf. the arguments in [BV1,
Proposition 1 and Lemma 3]).  We believe that, with the exception of $\A_5$,
the minimal dimension of an integer homology sphere on which the alternating
group $\A_{n+1}$ acts faithfully coincides with the bound $n-1$ in Theorem 1
(which is also the minimal dimension of a linear faithful action of $\A_{n+1}$
on a sphere); if so, this gives an independent proof of Theorem 1 for integer
homology spheres. This is indeed the case for various small values of $n$; we
note, however, that $\A_6$ acts faithfully on a mod 2 homology 3-sphere ([Zm1])
but not on an integer one (whereas e.g. for $\A_8$  the two minimal dimensions
coincide). Finite simple groups which admit actions on low-dimensional homology
spheres are considered in [MZ1],[MZ3] (dimension three) and [MZ2] (dimension
four), see also the survey [Zm2].

\bigskip

{\bf 2. Proof of Theorem 2}

\medskip

We consider a faithful action of a finite quotient $G$ of ${\rm SL}(n,\Bbb Z)$ on a
mod 2 homology $m$-sphere $M$, with $m < n-1$ and $n>2$, and have to show that $G$ is
trivial. We denote by $U$ the
kernel of the projection from ${\rm SL}(n,\Bbb Z)$ to $G$, with ${\rm SL}(n,\Bbb
Z)/U \cong G$.  We consider also the action of ${\rm SL}(n,\Bbb Z)$ on $M$
induced by the action of $G$, with $U$ acting trivially.  Since
${\rm SL}(n,\Bbb Z)$ is perfect for $n>2$ also $G$ is perfect (does not admit a
surjection onto a non-trivial abelian group), in particular $G$ and ${\rm SL}(n,\Bbb
Z)$ act orientation-preservingly on $M$. Arguing by contradiction, we will assume
in the following that $G$ is non-trivial and hence non-solvable.

\medskip

Since $U$ is a finite-index subgroup of ${\rm SL}(n,\Bbb Z)$ and $n > 2$, by
the congruence subgroup property (see [M] or [BMS])  $U$ contains a congruence
subgroup $C(k)$, for some positive integer $k$; so  $C(k)$ is the kernel of the
canonical homomorphism ${\rm SL}(n,\Bbb Z) \to {\rm SL}(n,\Z/k\Z)$ which  is
surjective (see e.g. [N, II.21]), and we have an exact sequence

$$1 \to  C(k)  \to  {\rm SL}(n,\Bbb Z) \to {\rm SL}(n,\Z/k\Z) \to 1.$$

In particular, there is a surjection from the finite group ${\rm SL}(n,\Bbb
Z)/C(k) \cong {\rm SL}(n,\Bbb Z/k\Bbb Z)$ to ${\rm SL}(n,\Bbb Z)/U \cong G$
which we denote by $\Phi: {\rm SL}(n,\Z/k\Z) \to G$.

\medskip

If $k = p_1^{r_1} \ldots p_s^{r_s}$ is the prime decomposition, then it is
well-known that
$${\rm SL}(n,\Bbb Z/k\Bbb Z) \; \cong \; {\rm SL}(n,\Bbb Z/p_1^{r_1}\Z)
\times \ldots \times {\rm SL}(n,\Bbb Z/p_s^{r_s}\Z)$$ (see [N, Theorem VII.11]
or [F, proof of Lemma 3.5.5.1]). Now the restriction of $\Phi: {\rm
SL}(n,\Z/k\Z) \to G$ to some factor
${\rm SL}(n,\Bbb Z/p_i^{r_i}\Z) = {\rm SL}(n,\Bbb Z/p^r\Z)$
has to be
non-trivial and induces a surjection $\Phi_0: {\rm SL}(n,\Z/p^r\Z) \to G_0$ onto a
perfect non-solvable subgroup $G_0$ of $G$; we denote by $U_0$ the kernel of $\Phi_0$
(the elements of ${\rm SL}(n,\Z/p^r\Z)$ acting trivially on $M$).

\medskip

Let $K$ denote the kernel of the canonical surjection ${\rm SL}(n,\Z/p^r\Z) \to
{\rm SL}(n,\Z/p\Z)$, so $K$ consists of all matrices in ${\rm SL}(n,\Z/p^r\Z)$
which are congruent to the identity matrix $I=I_n$ when entries are taken mod
$p$. By performing the binomial expansion of $(I + pA)^{p^{r-1}}$ one checks
that $K$ is a $p$-group, in particular $K$ is solvable (and the only
non-abelian factor in a composition series of  ${\rm SL}(n,\Z/p^r\Z)$ is the
simple group ${\rm PSL}(n,\Z/p\Z)$).

\medskip

Let $K_0$ denote the kernel of the surjection from ${\rm SL}(n,\Z/p^r\Z)$ to
the central quotient ${\rm PSL}(n,\Z/p\Z)$ of ${\rm SL}(n,\Z/p\Z)$; also $K_0$
is solvable and, since $n > 2$, ${\rm PSL}(n,\Z/p\Z)$ is a non-abelian simple
group. We will show that there exists some element $u$ in ${\rm
SL}(n,\Z/p^r\Z)$ which acts trivially on $M$ (i.e., $u \in U_0$) and projects
to a non-central element in ${\rm SL}(n,\Z/p\Z)$. Then $u$ projects
non-trivially also to the  central quotient ${\rm PSL}(n,\Z/p\Z)$ and hence the
normal subgroup $U_0$ of ${\rm SL}(n,\Z/p^r\Z)$ surjects onto the simple group
${\rm PSL}(n,\Z/p\Z)$. Considering the two exact sequences

$$ 1 \to K_0 \to {\rm SL}(n,\Z/p^r\Z) \to {\rm PSL}(n,\Z/p\Z) \to 1$$

$$ 1 \to U_0 \cap K_0 \to U_0 \to {\rm PSL}(n,\Z/p\Z) \to 1$$

and quotienting the first
by the second one concludes that  ${\rm SL}(n,\Z/p^r\Z)/U_0 \cong G_0$ is isomorphic
to the solvable group $K_0/U_0 \cap K_0$. This is a contradiction, and hence $G$ has
to be trivial.

\medskip

In order to find the element $u$, we distinguish the cases $p > 2$ and $p = 2$.

\bigskip

{\bf 2.1}  Suppose first that $p > 2$.

\medskip

Suppose that $n$ is odd. The subgroup $A$ of all diagonal matrices in ${\rm
SL}(n,\Z/p^r\Z)$ with all diagonal entries equal to $\pm 1$ is isomorphic to
$(\Z_2)^{n-1}$ and does not contain the central involution $-I$; note that this
subgroup injects under the canonical projection into ${\rm SL}(n,\Z/p\Z)$ and
hence also into its central quotient group ${\rm PSL}(n,\Z/p\Z)$ which is a
non-abelian simple group.  By general Smith fixed point theory the group
$(\Z_2)^{n-1}$ does not act faithfully and orientation-preservingly on a mod 2
homology sphere $M$ of dimension $m < n-1$ (see [P, Theorem 3.3]), so some
involution $u$ in $A$ acts trivially on $M$.

\medskip

If $n$ is even  then $A$ contains the central involution $-I$. If the central
involution acts non-trivially on $M$ then some other involution in $A$ has to
act trivially and we are done. There remains the case that the central
involution $-I$ acts trivially on $M$. In this case we appeal to [P] where it
is shown in the proof of Theorem 1.1 that there exists some non-central
involution represented by a diagonal matrix with diagonal entries $\pm1$ which
acts trivially on $M$.

\bigskip

{\bf 2.2}  Now suppose that $p = 2$.

\medskip

Denote by $E_{i,j}$ the $n \times n$-matrix with all entries zero except for
the $(i,j)$-entry which is equal to one. Let $A$ be  the subgroup of ${\rm
SL}(n,\Z/2^r\Z)$ generated by the elementary matrices  $I+E_{1,2}, \;
I+E_{1,3}, \ldots,\; I+E_{1,n}$. These matrices have order $2^r$ and commute,
so $A$ is isomorphic to $(\Z_{2^r})^{n-1}$. Again by Smith theory (see
[P, Theorem 3.3]), the subgroup  $(\Z_2)^{n-1}$ of $A$ does not act faithfully
on $M$. It follows easily that there is an element $u$ of order $2^r$ in $A$
which acts trivially on $M$. Under the canonical projection from ${\rm
SL}(n,\Z/2^r\Z)$ onto the simple group ${\rm SL}(n,\Z/2\Z) = {\rm
PSL}(n,\Z/2\Z)$, the element $u$ is mapped non-trivially.

\bigskip

This finishes the proof of Theorem 2.

\bigskip  \bigskip

Note (june 2008).  Parwani's result in [P] as well as the main result in a
preliminary version of the present paper (published in april 2006 in
arXiv:math) are formulated for continuous actions. Recently (march 2008) a
paper of Bridson and Vogtmann [BV2] appeared in arXiv:math, pointing out some
gap in Parwani's paper for the case of continuous actions (see Remarks 4.16 and
4.17 of the paper of Bridson and Vogtmann). Since Parwani's arguments remain
valid for smooth actions, we replaced continuous actions by smooth ones in the
present paper (and refer to the paper of Bridson and Vogtmann for the case of
continuous actions).

\bigskip
\bigskip

\centerline {\bf References}

\bigskip

\item {[BMS]} H.Bass, J.Milnor, J.P.Serre, {\it The congruence subgroup
property for $SL_n$ ($n \ge 3$) and  $SP_{2n}$ ($n \ge 2$).}    Inst.
Hautes Etudes Sci. Publ. Math. 33, 59-137 (1967)

\item {[BV1]} M.R.Bridson, K.Vogtmann, {\it  Homomorphisms from automorphism groups
of free groups.}  Bull. London Math. Soc. 35, 785-792  (2003)

\item {[BV2]} M.R.Bridson, K.Vogtmann, {\it  Actions of automorphism groups of
free groups on homology spheres and acyclic manifolds.}  arXiv:0803.2062

\item {[F]} B.Fine, {\it  Algebraic Theory of the Bianchi Groups.}
Monographs and Textbooks in Pure and Applied Mathematics, vol. 129,
Marcel Dekker, New York 1989

\item {[FM]} B.Farb, H.Masur, {\it Superrigidity and mapping class
groups.} Topology 37, 1169-1176 (1998)

\item {[FS]} B.Farb, P.Shalen, {\it Real-analytic actions of lattices.}
Invent. math. 135, 273-296  (1999)

\item {[M]} J.Mennicke, {\it Finite factor groups of the unimodular groups.}
Ann. Math. 81, 31-37 (1965)

\item {[MZ1]} M.Mecchia, B.Zimmermann, {\it On finite groups acting on
$\Z_2$-homology 3-spheres.} Math. Z. 248,  675-693  (2004)

\item {[MZ2]} M.Mecchia, B.Zimmermann, {\it On finite simple and
nonsolvable groups acting on homology 4-spheres.}
Topology Appl. 153, 2933-2942 (2006)

\item {[MZ3]} M.Mecchia, B.Zimmermann, {\it On finite simple groups acting
on integer and mod 2 homology 3-spheres.}  J. Algebra 298, 460-467 (2006)

\item {[N]} M.Newman, {\it Integral Matrices.} Pure and Applied Mathematics
Vol.45,  Academic Press 1972

\item {[P]} K.Parwani, {\it Actions of SL$(n,\Bbb Z)$ on homology
spheres.}  Geom. Ded. 112, 215-223 (2005)

\item {[We]} S.Weinberger, {\it SL$(n,\Bbb Z)$ cannot act on small tori.}
AMS/IP Studies in Advanced Mathematics Volume 2 (Part 1), 406-408 (1997)

\item {[Wi]} D.Witte, {\it Arithmetic groups of higher $\Bbb Q$-rank
cannot act on 1-manifolds.}  Proc. Amer. Math. Soc. 122, 333-340 (1994)

\item {[Z]} R.Zimmer, {\it Actions of semisimple groups and discrete subgroups.}
Proceedings of the I.C.M.,  Berkeley 1986,  1247-1258

\item {[Zm1]} B.Zimmermann, {\it Cyclic branched coverings and homology
3-spheres with large group actions.}  Fund. Math. 184, 343-353
(2004)

\item {[Zm2]} B.Zimmermann, {\it Some results and conjectures on finite
groups acting on homology spheres.} Sib. Electron. Math. Rep. 2,
233-238 (2005)   (http://semr.math.nsc.ru)

\bye